\documentclass[a4paper, oneside]{amsart}
\RequirePackage{amsmath}
\RequirePackage{bm}
\RequirePackage{amssymb}
\RequirePackage{upref}
\RequirePackage{amsthm}
\RequirePackage{enumerate}
\RequirePackage{pb-diagram}
\RequirePackage{amsfonts}
\RequirePackage[mathscr]{eucal}
\RequirePackage{verbatim}
\RequirePackage{xr}
\RequirePackage[pdftex]{graphicx}
\RequirePackage[pdftex]{floatflt}

\newcommand{\al}{\alpha}
\newcommand{\bet}{\beta}
\newcommand{\ga}{\gamma}
\newcommand{\de}{\delta }
\newcommand{\e}{\epsilon}

\newcommand{\f}{\varphi}
\newcommand{\h}{\eta}

\newcommand{\n}{\nu}

\newcommand{\s}{\sigma}
\newcommand{\x}{\xi}

\newcommand{\D}{\varDelta}
\newcommand{\F}{\varPhi}

\newcommand{\di}[1]{#1\nobreakdash-\hspace{0pt}dimensional}
\newcommand{\nbdd}{\nobreakdash--}

\newcommand{\fu}[3]{#1\hspace{0pt}_{|_{#2_#3}}}
\newcommand{\fv}[2]{#1\hspace{0pt}_{|_{#2}}}

\newcommand{\so}{{\mc S_0}}

\newcommand{\const}{\tup{const}}

\newcommand{\msp[1]}[1]{\mspace{#1mu}}

\newcommand{\R}[1][n+1]{{\protect\mathbb R}^{#1}}

\newcommand{\N}{{\protect\mathbb N}}

\newcommand{\eR}{\stackrel{\lower1ex \hbox{\rule{6.5pt}{0.5pt}}}{\msp[3]\R[]}}
\newcommand{\eN}{\stackrel{\lower1ex \hbox{\rule{6.5pt}{0.5pt}}}{\msp[1]\N}}
\newcommand{\eO}{\stackrel{\lower1ex
\hbox{\rule{6pt}{0.5pt}}}{\msc O}}

\DeclareMathOperator{\graph}{graph}

\newcommand\pde[2]{\frac {\partial#1}{\partial#2}}
 
\newcommand{\un}{\infty}

\newcommand{\set}[2]{\{\,#1\colon #2\,\}}
\newcommand{\uu}{\cup}

\newcommand{\uuu}{\bigcup}

\newcommand{\uud}{ \stackrel{\lower 1ex \hbox {.}}{\uu}}
\newcommand{\uuud}[1]{ \stackrel{\lower 1ex \hbox {.}}{\uuu_{#1}}}
\newcommand\su{\subset}

\newcommand{\sminus}[1][28]{\raise 0.#1ex\hbox{$\scriptstyle\setminus$}}

\newcommand\ti{\times }

\newcommand{\abs}[1]{\lvert#1\rvert}

\newcommand{\norm}[1]{\lVert#1\rVert}

\newcommand\ch[3]{\varGamma_{#1#2}^#3}
\newcommand\cha[3]{{\bar\varGamma}_{#1#2}^#3}
\newcommand\chc[3]{{\tilde\varGamma}_{#1#2}^#3}
\newcommand{\riem}[4]{R_{#1#2#3#4}}
\newcommand{\riema}[4]{{\bar R}_{#1#2#3#4}}

\newcommand{\tbf}{\textbf}
\newcommand{\tit}{\textit}

\newcommand{\tup}{\textup}

\newcommand{\mc}{\protect\mathcal}
\newcommand{\msc}{\protect\mathscr}

\providecommand{\bysame}{\makebox[3em]{\hrulefill}\thinspace}

\newcommand{\ci}{\cite}
\newcommand{\bib}{\bibitem}

\newcommand{\bt}{\begin{thm}}
\newcommand{\bl}{\begin{lem}}
\newcommand{\bc}{\begin{cor}}
\newcommand{\bd}{\begin{definition}}
\newcommand{\bpp}{\begin{prop}}
\newcommand{\br}{\begin{rem}}
\newcommand{\bn}{\begin{note}}
\newcommand{\be}{\begin{ex}}
\newcommand{\bes}{\begin{exs}}
\newcommand{\bb}{\begin{example}}
\newcommand{\bbs}{\begin{examples}}
\newcommand{\ba}{\begin{axiom}}

\newcommand{\et}{\end{thm}}
\newcommand{\el}{\end{lem}}
\newcommand{\ec}{\end{cor}}
\newcommand{\ed}{\end{definition}}
\newcommand{\epp}{\end{prop}}
\newcommand{\er}{\end{rem}}
\newcommand{\en}{\end{note}}
\newcommand{\ee}{\end{ex}}
\newcommand{\ees}{\end{exs}}
\newcommand{\eb}{\end{example}}
\newcommand{\ebs}{\end{examples}}
\newcommand{\ea}{\end{axiom}}

\newcommand{\bp}{\begin{proof}}
\newcommand{\ep}{\end{proof}}
\newcommand{\eps}{\renewcommand{\qed}{}\end{proof}}

\newcommand{\bal}{\begin{align}}

\newcommand{\bi}[1][1.]{\begin{enumerate}[\upshape #1]}
\newcommand{\bia}[1][(1)]{\begin{enumerate}[\upshape #1]}
\newcommand{\bin}[1][1]{\begin{enumerate}[\upshape\bfseries #1]}
\newcommand{\bir}[1][(i)]{\begin{enumerate}[\upshape #1]}
\newcommand{\bic}[1][(i)]{\begin{enumerate}[\upshape\hspace{2\cma}#1]}
\newcommand{\bis}[2][1.]{\begin{enumerate}[\upshape\hspace{#2\parindent}#1]}
\newcommand{\ei}{\end{enumerate}}

\newcommand\ndots{\raise 0.47ex \hbox {,}\hskip0.06em\cdots %
     \raise 0.47ex \hbox {,}\hskip0.06em} 

\newcommand{\q}{\quad}

\newcommand\nd{\noindent}

\newskip\Csmallskipamount                                                
\Csmallskipamount=\smallskipamount
\newskip\Cmedskipamount
\Cmedskipamount=\medskipamount
\newskip\Cbigskipamount
\Cbigskipamount=\bigskipamount

\newcommand\cvs{\vspace\Csmallskipamount}   
\newcommand\cvm{\vspace\Cmedskipamount}
\newcommand\cvb{\vspace\Cbigskipamount}

\newskip\csa
\csa=\smallskipamount

\newskip\cma
\cma=\medskipamount

\newskip\cba
\cba=\bigskipamount

\newdimen\spt
\newcommand\citem{\cvs\advance\itemno by
1{(\romannumeral\the\itemno})\hskip3pt}
\newcommand{\bitem}{\cvm\nd\advance\itemno by
1{\bf\the\itemno}\hspace{\cma}}

\newcount\itemno
\newcommand{\las}[1]{\label{S:#1}}

\newcommand{\lae}[1]{\label{E:#1}}

\newcommand{\rs}[1]{Section~\ref{S:#1}}

\newcommand{\re}[1]{\eqref{E:#1}}

\RequirePackage{amsmath}
\RequirePackage{amssymb}

\RequirePackage{upref}
\RequirePackage{amsthm}
\usepackage{amsfonts}

\RequirePackage{enumerate}
\usepackage[mathscr]{eucal}

\theoremstyle{plain}
\newtheorem{thm}{Theorem}[section]
\newtheorem{lem}[thm]{Lemma}
\newtheorem{prop}[thm]{Proposition}
\newtheorem{cor}[thm]{Corollary}

\theoremstyle{definition}
\newtheorem{rem}[thm]{Remark}
\newtheorem{definition}[thm]{Definition}
\newtheorem{example}[thm]{Example}
\newtheorem{ex}[thm]{Exercise}

\swapnumbers
\theoremstyle{remark}

\numberwithin{equation}{section}

\usepackage{hyperref}

\begin{document}
\title[Constant mean curvature foliation of space-time]{On the foliation of
space-time by constant mean curvature hypersurfaces}

\author{Claus Gerhardt}
\address{Ruprecht-Karls-Universit\"at, Institut f\"ur Angewandte Mathematik,
Im Neuenheimer Feld 294, 69120 Heidelberg, Germany}
\email{gerhardt@math.uni-heidelberg.de}

%
\date{March 21, 1999}
%


\begin{abstract} We prove that the mean curvature $\tau$ of the slices given by a
constant mean curvature foliation can be used as a time function, i.e. $\tau$ is
smooth with non-vanishing gradient.
\end{abstract}
\maketitle
\tableofcontents
\setcounter{section}{-1}
\section{Introduction} 

\cvb
In \ci{cg1} it is proved that a globally hyperbolic Lorentzian manifold $N$ with a
compact Cauchy hypersurface can be foliated by constant mean curvature
hypersurfaces if the \tit{big bang} and \tit{ big crunch} hypotheses are valid, and
if the \tit{time-like convergence} condition holds, i.e. if
\begin{equation}\lae{0.1}
\bar R_{ \al\bet}\h^ \al\h^\bet\ge 0
\end{equation}
for all time-like vectorfields $(\h^ \al)$.

If we assume for simplicity that the mean curvature of the barriers provided by
the big bang resp. big crunch hypotheses tend to $-\un$ resp. $+\un$, the foliation
can be described as follows: Indicate by $M_\tau$ a closed hypersurface of mean
curvature $\tau$, then, the foliation consists of the uniquely determined $M_\tau$,
$0\ne \tau\in \R[]$, and of the set $\msc C_0$ of \tit{maximal} slices. If there is
more than one maximal slice, then $\msc C_0$ comprises a whole continuum of
maximal slices, which are all totally geodesic and the ambient metric is static in
$\msc C_0$.

The mean curvature of the slices of the foliation can be looked at as a function
$\tau$ on $N$, which is continuous as one easily checks. However, in order to use
$\tau$ as a new time function, $\tau$ has to be smooth with non-vanishing
gradient.

We prove that this is indeed the case in $\{\tau\ne 0\}$ and also globally, if there
is a maximal slice that is not totally geodesic or if the strict inequality is valid in
\re{0.1}. Evidently, if we have two (and more) maximal slices, then, $D\tau$
vanishes in the interior of
$\msc C_0$, and we shall also give an example of a foliation with exactly one
totally geodesic slice, where
$D\tau$ vanishes on that slice, i.e. the assumptions guaranteeing $\tau$ to be
smooth and $D\tau\ne 0$ are also sharp.

The paper is organized as follows: In  \rs{1} we introduce the notations
and common definitions we rely on.

The main theorem is proved in \rs{2}, while the counterexample is given in \rs{3}.

\cvb
\section{Notations and definitions}\las 1

\cvb

The main objective of this section is to state the equations of Gau{\ss}, Codazzi,
and Weingarten for space-like hypersurfaces $M$ in a \di {(n+1)} Lorentzian
manifold
$N$.  Geometric quantities in $N$ will be denoted by
$(\bar g_{ \al \bet}),(\riema  \al \bet \ga \de)$, etc., and those in $M$ by $(g_{ij}), (\riem
ijk)$, etc. Greek indices range from $0$ to $n$ and Latin from $1$ to $n$; the
summation convention is always used. Generic coordinate systems in $N$ resp.
$M$ will be denoted by $(x^ \al)$ resp. $(\x^i)$. Covariant differentiation will
simply be indicated by indices, only in case of possible ambiguity they will be
preceded by a semicolon, i.e. for a function $u$ in $N$, $(u_ \al)$ will be the
gradient and
$(u_{ \al \bet})$ the Hessian, but e.g., the covariant derivative of the curvature
tensor will be abbreviated by $\riema  \al \bet \ga{ \de;\e}$. We also point out that
\begin{equation}
\riema  \al \bet \ga{ \de;i}=\riema  \al \bet \ga{ \de;\e}x_i^\e
\end{equation}
with obvious generalizations to other quantities.

Let $M$ be a \tit{space-like} hypersurface, i.e. the induced metric is Riemannian,
with a differentiable normal $\n$ which is time-like.

In local coordinates, $(x^ \al)$ and $(\x^i)$, the geometric quantities of the
space-like hypersurface $M$ are connected through the following equations
\begin{equation}\lae{1.2}
x_{ij}^ \al= h_{ij}\n^ \al
\end{equation}
the so-called \tit{Gau{\ss} formula}. Here, and also in the sequel, a covariant
derivative is always a \tit{full} tensor, i.e.

\begin{equation}
x_{ij}^ \al=x_{,ij}^ \al-\ch ijk x_k^ \al+ \cha  \bet \ga \al x_i^ \bet x_j^ \ga.
\end{equation}
The comma indicates ordinary partial derivatives.

In this implicit definition the \tit{second fundamental form} $(h_{ij})$ is taken
with respect to $\n$.

The second equation is the \tit{Weingarten equation}
\begin{equation}
\n_i^ \al=h_i^k x_k^ \al,
\end{equation}
where we remember that $\n_i^ \al$ is a full tensor.

Finally, we have the \tit{Codazzi equation}
\begin{equation}
h_{ij;k}-h_{ik;j}=\riema \al \bet \ga \de\n^ \al x_i^ \bet x_j^ \ga x_k^ \de
\end{equation}
and the \tit{Gau{\ss} equation}
\begin{equation}
\riem ijkl=- \{h_{ik}h_{jl}-h_{il}h_{jk}\} + \riema  \al \bet\ga \de x_i^ \al x_j^ \bet
x_k^ \ga x_l^ \de.
\end{equation}

Now, let us assume that $N$ is a globally hyperbolic Lorentzian manifold with a
\tit{compact} Cauchy surface. 
$N$ is then a topological product $\R[]\times \mc S_0$, where $\mc S_0$ is a
compact Riemannian manifold, and there exists a Gaussian coordinate system
$(x^ \al)$, such that the metric in $N$ has the form 
\begin{equation}\lae{1.7}
d\bar s_N^2=e^{2\psi}\{-{dx^0}^2+\s_{ij}(x^0,x)dx^idx^j\},
\end{equation}
where $\s_{ij}$ is a Riemannian metric, $\psi$ a function on $N$, and $x$ an
abbreviation for the space-like components $(x^i)$, see \ci{GR},
\ci[p.~212]{HE}, \ci[p.~252]{GRH}, and \ci[Section~6]{cg1}.
We also assume that
the coordinate system is \tit{future oriented}, i.e. the time coordinate $x^0$
increases on future directed curves. Hence, the \tit{contravariant} time-like
vector $(\x^ \al)=(1,0,\dotsc,0)$ is future directed as is its \tit{covariant} version
$(\x_ \al)=e^{2\psi}(-1,0,\dotsc,0)$.

Let $M=\graph \fv u\so$ be a space-like hypersurface
\begin{equation}
M=\set{(x^0,x)}{x^0=u(x),\,x\in\mc S_0},
\end{equation}
then the induced metric has the form
\begin{equation}
g_{ij}=e^{2\psi}\{-u_iu_j+\s_{ij}\}
\end{equation}
where $\s_{ij}$ is evaluated at $(u,x)$, and its inverse $(g^{ij})=(g_{ij})^{-1}$ can
be expressed as
\begin{equation}\lae{2.10}
g^{ij}=e^{-2\psi}\{\s^{ij}+\frac{u^i}{v}\frac{u^j}{v}\},
\end{equation}
where $(\s^{ij})=(\s_{ij})^{-1}$ and
\begin{equation}\lae{2.11}
\begin{aligned}
u^i&=\s^{ij}u_j\\
v^2&=1-\s^{ij}u_iu_j\equiv 1-\abs{Du}^2.
\end{aligned}
\end{equation}
Hence, $\graph u$ is space-like if and only if $\abs{Du}<1$.

The covariant form of a normal vector of a graph looks like
\begin{equation}
(\n_ \al)=\pm v^{-1}e^{\psi}(1, -u_i).
\end{equation}
and the contravariant version is
\begin{equation}
(\n^ \al)=\mp v^{-1}e^{-\psi}(1, u^i).
\end{equation}
Thus, we have
\br Let $M$ be space-like graph in a future oriented coordinate system. Then, the
contravariant future directed normal vector has the form
\begin{equation}
(\n^ \al)=v^{-1}e^{-\psi}(1, u^i)
\end{equation}
and the past directed
\begin{equation}\lae{2.15}
(\n^ \al)=-v^{-1}e^{-\psi}(1, u^i).
\end{equation}
\er

In the Gau{\ss} formula \re{1.2} we are free to choose the future or past directed
normal, but we stipulate that we always use the past directed normal for reasons
that we have explained in \ci[Section 2]{cg5}.

Look at the component $ \al=0$ in \re{1.2} and obtain in view of \re{2.15}

\begin{equation}\lae{1.16}
e^{-\psi}v^{-1}h_{ij}=-u_{ij}- \cha 000\mspace{1mu}u_iu_j- \cha 0j0
\mspace{1mu}u_i- \cha 0i0\mspace{1mu}u_j- \cha ij0.
\end{equation}
Here, the covariant derivatives are taken with respect to the induced metric of
$M$, and
\begin{equation}
-\cha ij0=e^{-\psi}\bar h_{ij},
\end{equation}
where $(\bar h_{ij})$ is the second fundamental form of the hypersurfaces
$\{x^0=\const\}$.

An easy calculation shows
\begin{equation}
\bar h_{ij}e^{-\psi}=-\tfrac{1}{2}\dot\s_{ij} -\dot\psi\s_{ij},
\end{equation}
where the dot indicates differentiation with respect to $x^0$.

\cvb
\section{$\tau$ is a time function}\las{2}

\cvb
Let $M_0=M_{\tau_0}$ be a hypersurface of constant mean curvature $\tau_0$, and
let
$(x^ \al)$ be a future oriented, normal Gaussian coordinate system relative to
$M_0$, i.e.
\begin{equation}
d\bar s_N^2=-{dx^0}^2+\s_{ij}(x^0,x)dx^idx^j
\end{equation}
in a tubular neighbourhood $\mc U$ of $M_0$; here, $x$ is an abbreviation for the
spatial components $(x^i)$.

Since we have a continuous foliation of $N$, the constant mean curvature slices
contained in $\mc U$ can be written as graphs over $M_0$
\begin{equation}
M_\tau=\set {(x^0,x^i)}{x^0=u(\tau,x^i)}
\end{equation}
with a continuous function $u\in C^0(I\ti M_0)$, $I$ an open interval containing
$\tau_0$.

Let us assume for the moment that $u$ is smooth---we shall see later that this is
indeed the case, at least for $\tau_0\ne 0$---, and let us define the
transformation
\begin{equation}
\varPhi(\tau,x^i)=(u(\tau,x^i),x^i).
\end{equation}
Then, we have
\begin{equation}
\det D\F=\pde u\tau\equiv \dot u.
\end{equation}

In view of the monotonicity of the foliation, cf. \ci[Lemma 7.2]{cg1}, we know that
$\dot u\ge 0$. Thus, if we could show that $u$ is smooth and $\dot u$ strictly
positive, we would obtain that $\F$ is a diffeomorphism, and hence, that $\tau$ is
smooth with non-vanishing gradient.

Let us first show that $u$ is smooth.

\bl
Let $\tau_0\in \R[]$ be such  that $M_{\tau_0}$ is not totally geodesic or assume
that the strict inequality is valid in \re{0.1}. Consider a tubular neighbourhood
$\mc U$ of
$M_0=M_{\tau_0}$, and choose normal Gaussian coordinates as above, such that
the
$M_\tau\su\mc U$ are graphs
\begin{equation}
M_\tau=\graph \fu {u(\tau)}M0.
\end{equation}
Then, $u$ is smooth in $I\ti M_0$, where $I=(\tau_0-\e,  \tau_0+\e)$ for small
$\e=\e( \tau_0)>0$.
\el

\bp
Let $u_0=u(\tau_0)$, which is incidentally identical zero due to our choice of
coordinates,  and define the operator
$G$
\begin{equation}
G(\tau,\f)=H(\f)-\tau,
\end{equation}
where $H(\f)$ is an abbreviation for the mean curvature of $\graph \fu {\f}M0$.

It is well-known, cf. \ci[Section 5]{cb}, that $D_2G(\tau_0,u_0)\f$ equals
\begin{equation}
-\D\f+\f\{\norm A^2+\bar R_{ \al \bet}\n^ \al\n^ \bet\},
\end{equation}
where $\norm A^2=h_{ij}h^{ij}$, and the Laplacian, the second fundamental form
and $\n$ are evaluated with respect to $M_0$.

Hence, $D_2G$ is invertible at $(\tau_0,u_0)$ in view of our assumptions,  and thus,
the
\tit{implicit function theorem} is applicable to yield the existence of a smooth
function
$u=u(\tau,x)$ satisfying
\begin{equation}
G(\tau,u(\tau,x))=0
\end{equation}
for $\abs{\tau-\tau_0}<\e$, $\e$ small; of course,  the new smooth function
coincides with the original continuous function.
\ep

It remains to demonstrate that $\dot u$ is strictly positive.

We use the same setting as above, i.e. we work in a tubular
neighbourhood $\mc U$ of a slice $M_{\tau_0}$ with corresponding normal
Gaussian coordinates, and look at the $M_\tau\su \mc U$, which are then graphs
over
$M_{\tau_0}$.

\bl
Let $M_\tau=\graph u(\tau)$ and $M_{\bar \tau}=\graph  u(\bar \tau)$, then, there
exists a constant $c=c(\mc U)$ such that
\begin{equation}\lae{2.9}
\abs{\tau-\bar \tau}\le c\inf_{M_{\tau_0}}\abs{u-\bar u},
\end{equation}
where $\bar u=u(\bar \tau)$.
\el

\bp
Assume for simplicity $\bar \tau<\tau$. We want to employ the relation
\re{1.16}---now with $\psi\equiv 0$---, but the covariant derivatives of $u$
should be expressed with respect to the metric $(\s_{ij}(u,x))$; recall that
\begin{equation}
g_{ij}=-u_iu_j+\s(u,x).
\end{equation}

An easy calculation reveals that
\begin{equation}
u_{ij}=v^{-2}u_{;ij}\msp,
\end{equation}
where the semicolon indicates the covariant derivatives with respect to the
metric
$(\s_{ij}(u,x))$. Let $\chc ijk$ be the corresponding Christoffel symbols and
indicate by a comma ordinary partial derivatives, then, we obtain
\begin{equation}\lae{2.22}
h_{ij}=-v^{-1}u_{,ij}+v^{-1}\chc ijk u_k-v\cha ij0
\end{equation}
with an analogous relation valid for $\bar u$; note that the other Christoffel
symbols in \re{1.16} vanish, since $\psi\equiv 0$.

Let $x_0\in M_{\tau_0}$ be such that
\begin{equation}
(u-\bar u)(x_0)=\inf_{M_{\tau_0}}(u-\bar u).
\end{equation}

We now observe that in $x_0$
\begin{equation}
Du=D\bar u\q \tup{and}\q D^2u\ge D^2\bar u,
\end{equation}
and deduce from \re{2.22}, and the corresponding relation for $M_{\bar \tau}$,
\begin{equation}
\tau-\bar \tau\le f(x_0,u(x_0))-f(x_0,\bar u(x_0))
\end{equation}
with a smooth function $f$, i.e.
\begin{equation}
\tau-\bar \tau\le c[u(x_0)-\bar u(x_0)].
\end{equation}
\ep

\br
An estimate of the form \re{2.9} is valid in \tit{any} Gaussian coordinate system
in which the slices of the foliation are represented as graphs. The constant $c$
depends on the compact set containing the slices under consideration.
\er

Thus, we have proved

\bt
The function $\tau$ is smooth with non-vanishing gradient in $\{\tau\ne 0\}$. If
there is a maximal slice that is not totally geodesic or if the strict inequality is
valid in \re{0.1}, then
$\tau$ is a globally defined time function.
\et

We note that by assumption there are maximal slices in the foliation.

\cvb
\section{A counterexample}\las{3}

\cvb
Let $S^n$ be the standard unit sphere with metric $\s_{ij}$, and for $t\in
I=(-\e,\e)$ let
\begin{equation}
f=f(t)=-\int_0^t\frac{s^3}{\e^2-s^2}\,ds.
\end{equation}
Then, we define $N=I\ti S^n$ with metric
\begin{equation}
d\bar s_N^2=-dt^2+e^{2f}\s_{ij}\msp dx^i dx^j,
\end{equation}
where the coordinates $(t,x^i)$ are supposed to be future oriented.

The level hypersurfaces $\{t=\const\}$ are all totally umbilical, and their mean
curvature is equal to
\begin{equation}
\tau=H=-n\dot f=n\frac{t^3}{\e^2-t^2}\raise 2pt \hbox{.}
\end{equation}

Thus, we see that the big bang and big crunch hypotheses are satisfied and that
$D\tau$ vanishes on  the unique totally geodesic maximal slice.

It remains to verify  the time-like convergence condition.

We first note, that for $\e\le 1$
\begin{equation}\lae{3.4}
\ddot f+\dot
f^2=-\frac{3t^2}{\e^2-t^2}-\frac{2t^4}{(\e^2-t^2)^2}+\frac{t^6}{(\e^2-t^2)^2}\le
0,
\end{equation}
which will be the core inequality to estimate the Ricci tensor.

Secondly, let us write the metric as a conformal metric
\begin{equation}
\bar g_{ \al \bet}=e^{2\psi}g_{ \al \bet}
\end{equation}
in coordinates $(x^0,x^i)$, where

\begin{align}
g_{ \al \bet}dx^ \al dx^ \bet&=-{dx^0}^2+\s_{ij}\msp dx^idx^j,\\
\frac{dx^0}{dt}&=e^{-f(t)},
\end{align}
and $\psi(x^0)=f(t)$.

The metric $(g_{ \al \bet})$ is a product metric and, therefore, the only non-zero
components of its Ricci tensor $(R_{ \al \bet})$ are of the form $R_{ij}$ and coincide
with the components of the Ricci tensor of $S^n$.

The Ricci tensor $(\bar R_{ \al \bet})$ of the metric $(\bar g_{ \al \bet})$ is connected
to $(R_{ \al \bet})$ through the relation
\begin{equation}
\begin{split}
\bar R_{ \al \bet}&=R_{ \al \bet}-(n-1)[\psi_{ \al \bet}-\psi_ \al\psi_ \bet]\\
&\q\,-g_{ \al \bet}[\D\psi+(n-1)\abs
{D\psi}^2],
\end{split}
\end{equation}
where covariant derivatives and the norm are calculated with respect to
$(g_{ \al \bet})$; note that $\abs{D\psi}^2$ can be negative.

Hence, we conclude
\begin{equation}
\begin{aligned}
\bar R_{00}&=-n\ddot\psi =-ne^{2f}[\ddot f+\dot f^2],\\
\bar R_{0i}&=0,
\end{aligned}
\end{equation}
and
\begin{equation}
\begin{aligned}
\bar R_{ij}&=R_{ij}+\s_{ij}[\ddot\psi+(n-1)\dot\psi^2],\\
&=R_{ij}+\s_{ij}[\ddot f+n\dot f^2]e^{2f}.
\end{aligned}
\end{equation}

Let $(\h^ \al)$ be a time-like vector with component $\h^0=1$, so that
\begin{equation}
\s_{ij}\h^i\h^j<1.
\end{equation}

Then, taking into account that $\ddot f\le 0$, we derive

\begin{equation}
\begin{split}
\bar R_{ \al \bet}\h^ \al\h^ \bet&=-ne^{2f}[\ddot f+\dot f^2]+R_{ij}\h^i\h^j\\
&\q\,+\s_{ij}\h^i\h^j[\ddot f+n\dot f^2]e^{2f}\\
&\ge ne^{2f}[\ddot f+\dot f^2][-1 +\s_{ij}\h^i\h^j],
\end{split}
\end{equation}
and we conclude that the time-like convergence condition is satisfied in view of
\re{3.4}.

\end{document}